\documentclass[review]{elsarticle}

\usepackage{lineno,hyperref}
\usepackage{amsmath}
\usepackage{mymacros}
\DeclareMathOperator*{\starspan}{span}

\usepackage{algorithm}
\usepackage{algorithmic} 
\usepackage{enumitem}
\newlist{steps}{enumerate}{1}
\setlist[steps, 1]{label = Step \arabic*:}

\newtheorem{theorem}{Theorem}
\newtheorem{lemma}{Lemma}

\usepackage{tikz}
\usepackage{graphicx}
\usepackage{rotating}
\usepackage{pgfplots}
\usepackage{caption}
\usepackage{subcaption}
\modulolinenumbers[5]

\journal{journal}

\begin{document}

\begin{frontmatter}

\title{Multivariate moment matching for model order reduction of quadratic-bilinear systems using error bounds}

\author[mymainaddress]{Muhammad Altaf Khattak\fnref{fn1}}
\author[mymainaddress]{Mian Ilyas Ahmad\corref{mycorrespondingauthor}\fnref{fn1}}
\ead{m.ilyas@rcms.nust.edu.pk}

\author[mysecondaryaddress]{Lihong Feng}
\author[mysecondaryaddress]{Peter Benner}

\fntext[fn1]{This author is supported by \href{https://hec.gov.pk/english/pages/home.aspx}{HEC} Pakistan under NRPU Project ID 10176.}

\cortext[mycorrespondingauthor]{Corresponding author}
\address[mymainaddress]{Research Center for Modelling and Simulation, NUST H-12, Islamabad 44000 Pakistan}
\address[mysecondaryaddress]{Computational Methods for Systems and Control, Max Planck Institute Magdeburg, Sandtorstrasse~1, 39106 Germany}

\begin{abstract}
We propose an adaptive moment-matching framework for model order reduction of quadratic-bilinear descriptor systems. In this framework, an important issue is the selection of those shift frequencies where moment-matching is to be achieved.  Often, the choice is random or linked to the linear part of the nonlinear system. In this paper, we extend the use of an existing a posteriori error bound for general linear time invariant systems to quadratic-bilinear systems and develop a greedy-type framework to select a good choice of interpolation points for the construction of the projection matrices. The results are compared with standard quadratic-bilinear projection methods and we observe that the approximations obtained by the proposed method yield high accuracy.
\end{abstract}

\end{frontmatter}

\newpage

\section{Introduction}

There are different applications where the dynamics of the system can be represented by quadratic-bilinear differential algebraic equations (QBDAEs). These include simulation of distribution networks \cite{Grundel2013}, fluid flow problems \cite{morKunV08} and nonlinear VLSI circuits \cite{morPhi03, morGu11}. In addition, a large class of nonlinear systems can be written in quadratic-bilinear form by using exact transformations \cite{morGu11}. Most of these applications involve large number of equations, which make simulation, control and optimization computationally inefficient. A remedy to this issue is the use of model order reduction (MOR).

We consider the problem of MOR for a single-input single-output quadratic-bilinear descriptor system of the form:
\begin{equation} \label{original_qbdae}
\begin{aligned}
E\dot{x}(t)&=Ax(t)+Nx(t)u(t)+Qx(t)\otimes x(t)+Bu(t),\\
~~\quad y(t) &=Cx(t),
\end{aligned}
\end{equation}

where $E, ~\!A, ~\!N\in \R^{n\times n}$, $Q\in \R^{n\times n^2}$, $B,~\!C^T\in \R^{n}$ are the coefficient matrices and vectors. $x(t)\in \R^n$ is the state vector and $u(t),~\!y(t)\in \R$ are the input and output of the system. The matrix $E$ may or may not be singular but the pencil is assumed to be regular, i.e., $\lambda E-A$ is singular only for finitely many values of $\lambda\in \C$ \cite{morFre11}.

The goal of MOR is to construct a reduced system of dimension $r\ll n$:

\begin{equation}\label{reduced_qbdae}
\begin{aligned}
E_r\dot{x}_r(t)&=A_r x_r(t)+N_r x_r(t) u(t)+Q_r x_r(t)\otimes x_r(t)+B_r u(t),\\
~~\quad y_r(t) &=C_r x_r(t),
\end{aligned}
\end{equation}

with the output response $y_r(t)$ approximately equal to $y(t)$. In case of linear systems (where $Q$ and $N$ are null matrices), there are various techniques in the literature to compute reduced-order models (ROMs), cf., \cite{morAnt05, morBauBF14}. Among these methods, projection-based moment-matching methods \cite{morGri97,morAntSG01} are well used and are recently extended to quadratic-bilinear systems \cite{morGu11, morBenB15,morAhmBJ16}. Projection involves approximating the state vector $x(t)$ in an $r$-dimensional subspace spanned by the column vectors of $V\in \Rnr$, so that the residual in the state equation is orthogonal to another $r$-dimensional subspace spanned by the column vectors of $W\in \Rnr$. That is, we approximate $x(t)\approx Vx_r(t)$ such that the Petrov-Galerkin orthogonality condition holds:

\begin{equation}
\begin{aligned}
\!\!\!&W^T\bigg(\!EV\dot{x}_r(t)\!-\!\big(AVx_r(t)+NVx_r(t)u(t)+QVx_r(t)\otimes Vx_r(t)+Bu(t)\big)\bigg)=0,\\
\!\!\!&\hat{y}(t) =CVx_r(t).
\end{aligned}
\end{equation}
If $W=V$, the projection is orthogonal and is often called one-sided projection, otherwise it is oblique and is called two-sided projection. The oblique projection framework leads to a set of reduced system matrices of the form:
\begin{equation}\label{reduced_mats}
\begin{aligned}
E_r=W^TEV,~~~A_r=&W^TAV,~~~Q_r=W^TQ(V\otimes V),~~~N_r=W^TNV,\\
&B_r =W^TB,~~~C_r = CV.
\end{aligned}
\end{equation}
In case of linear systems, a suitable choice of the basis matrices $V$ and $W$, implicitly ensure moment-matching, where moments are the coefficients of the series expansion of the transfer function at some predefined shift frequencies. Thus for projection-based moment-matching, the choice of $V$ and $W$ is related to the transfer function of the system. However, nonlinear systems have no universal input-output representation though for some classes of nonlinear systems, including the QBDAE system, it is possible to generalise the transfer function concept by utilising the Volterra theory \cite{nonlinearRugh81}, where the input-output relationship is represented by a set of high-order transfer functions. This makes the concept of moment-matching slightly complex in the nonlinear case, since the structure of the basis matrices $V$ and $W$ in \eqref{reduced_mats} now depends on multiple high-order transfer functions. To achieve moment-matching, some simplifications are made in the literature \cite{morGu11, morBenB15}  for computing the ROMs. For example, \cite{morBenB15} constructs $V$ and $W$ such that the reduced system matches the moments of the first- and second-order transfer functions. In \cite{morAhmBJ16}, simplified forms of high-order transfer functions are derived, which also enable the projection based techniques to match moments of high-order transfer functions. In addition, all the existing moment-matching/interpolation approaches \cite{morGu11, morBenB15, morAhmBJ16} are based on the simplification that the interpolation points for each frequency variable is the same. We discuss these results further in Section~\ref{section:qbmor}. 

Recently a new framework \cite{morBenGG16} for quadratic-bilinear systems has been proposed that is based on generalized Sylvester-type matrix equations. The approach involves truncated solution of two complex matrix equations to identify a good choice for the basis matrices $V$ and $W$. Another approach is the extension of the Loewner framework from linear/bilinear systems \cite{morMayA07,morIonA14} to quadratic bilinear systems \cite{morGosA15} . Also an indirect approach for MOR of the QBDAE system is proposed in \cite{morAhmFB15}, where the basis matrices are constructed from the bilinear part of the quadratic-bilinear system. In \cite{ahmad2018interpolatory}, the linear-bilinear part of the system is viewed as a linear parametric system and a posteriori error bound is used to select the interpolation points and construct the basis matrices adaptively. All these techniques are using the first two or three high-order transfer functions  and their structure is different from the one identified in \cite{morBenB15}. Since our target is moment-matching for QBDAEs, we will mainly focus on the two-sided moment-matching technique of Benner and Breiten \cite{morBenB15}.

In this paper, we identify a good choice of interpolation points for the quadratic-bilinear system by utilizing a greedy type framework based on error bounds for quadratic-bilinear systems motivated by the recently proposed error bound for linear parametric systems in \cite{morFenAB15}. Here we relax the restriction of using the same interpolation points for different frequency variables. The approach starts from some initial interpolation points that are iteratively updated to identify a set of interpolation points corresponding to the maximal values of certain error bounds. For each choice of interpolation points, we interpolate, not only, the original transfer function and its first derivative but also higher derivatives, so that the quadratic-bilinear system is well approximated. The iterations stops when the approximation error is less than the prescribed tolerance level. Each iteration contributes to constructing a better set of basis matrices $V$ and $W$, until a given error tolerance is achieved. The main difference from the work in \cite{ahmad2018interpolatory} is that the quadratic part of the system is also involved in basis construction in the proposed framework based on a posteriori error bound for quadratic-bilinear systems, whereas only the bilinear part is considered for the basis matrix computation in \cite{ahmad2018interpolatory}. The error estimator used in \cite{ahmad2018interpolatory} only estimates the error of the linear-bilinear part.

The remaining part of the paper is organized as follows. Section~2 reviews the existing projection based moment-matching techniques for quadratic-bilinear systems. Section~3 presents the error bound expressions for quadratic bilinear systems and Section~4 utilises these error bounds in a greedy-type algorithmto select interpolation points. Finally in Section~5, numerical results are shown for some benchmark examples.

\section{Background} \label{section:qbmor}
In this section, we briefly review the concept of moment-matching discussed in \cite{morBenB15, morAhmBJ16} for quadratic-bilinear systems. Before going into the details of nonlinear moment-matching, we begin with the structure of high-order transfer functions.

\subsection{Multivariate Transfer Functions}
The input-output representation for single input quadratic-bilinear systems can be expressed by the Volterra series expansion of the output $y(t)$ with quantities analogous to the standard convolution operator. That is, 
\be\label{io_time}
\begin{aligned}
y(t) = \sum_{k=1}^{\infty}\int_{0}^{t}\!\!\int_{0}^{t_1}\!\!\!\!\cdots\int_{0}^{t_{k-1}}h_k(t_1,\ldots,t_k)u(t-t_1)\cdots u(t-t_k)dt_k\cdots dt_1,
\end{aligned}
\ee
where it is assumed that the input signal is one-sided, i.e., $u(t)=0$ for $t<0$. In addition, each of the generalized impulse responses, $h_k(t_1,\ldots,t_k)$, also called the $k$-dimensional kernel of the subsystem, is assumed to be one-sided. 
In terms of the multivariate Laplace transform, the $k$-dimensional subsystem can be represented as,
\be \label{io_subsystem} 
Y_k(s_1,\ldots,s_k) = H_k(s_1,\ldots,s_k)U(s_1)\cdots U(s_k), 
\ee
where $H_k(s_1,\ldots,s_k)$ is the multivariate transfer function of the $k$-dimensional subsystem. The generalized transfer functions in the output expression \eqref{io_subsystem} are in the so-called triangular form \cite{nonlinearRugh81}. We denote the $k$-dimensional triangular form by $H_{tri}^{[k]}(s_1,\ldots,s_k)$. There are some other useful forms such as the symmetric form and the regular form of the multivariate transfer functions as discussed in \cite{nonlinearRugh81}. The triangular form is related to the symmetric form by the following expression
\be\label{tri2sym}
H_{sym}^{[k]}(s_1,\ldots,s_k)=\frac{1}{n!}\sum_{\pi(\cdot)}H_{tri}^{[k]}(s_{\pi(1)},\ldots,s_{\pi(k)}),
\ee
where the summation includes all $k!$ permutations of $s_1,\ldots,s_k$. Also, the triangular form can be connected to the regular form of the transfer function by using
\be\label{reg2tri}
H_{tri}^{[k]}(s_1,\ldots,s_k)=H_{reg}^{[k]}(s_1,s_1+s_2,\ldots,s_1+s_2+\cdots+s_k).
\ee
According to \cite{nonlinearRugh81}, the structure of the generalized symmetric transfer functions can be identified by the growing exponential approach.
The structure of these symmetric transfer functions for the first two subsystems of the quadratic-bilinear system \eqref{original_qbdae} can be written as
\be \label{sym_3}
\begin{aligned}
H_1(s_1)&=C(s_1E-A)^{-1}B,\\
H_2(s_1,s_2)&=C((s_1+s_2)E-A)^{-1}B(s_1,s_2),
\end{aligned}
\ee
here 
\be 
\begin{aligned}
&B(s_1,s_2) = :Q(x_1(s_1)\otimes x_1(s_2))+\frac{1}{2!}N(x_1(s_1)+x_1(s_2)),
\end{aligned}
\ee
in which $x_1(s):=(sE-A)^{-1}B$  and $Q$ satisfies $Q(u\otimes v)=Q(v\otimes u)$ for all $u,v\in \R^n$. Defining $x_2(s_1,s_2):= ((s_1+s_2)E-A)^{-1}B(s_1,s_2)$, the first two (first- and second-order) symmetric transfer functions can be written as
\be
\begin{aligned}
H_1(s_1)&=Cx_1(s_1),\\
H_2(s_1,s_2)&=Cx_2(s_1,s_2).
\end{aligned}
\ee

Before going into the partial differentiation of these multivariate transfer functions, we introduce the concept of matricization. The process of reshaping a tensor into a matrix is called matricization. In \cite{morBenB15},  the matrix $Q\in \R^{n\times n^2}$ is considered as the mode-1 matricization of a 3 dimensional tensor $\mathcal{Q}\in \R^{n\times n\times n}$. The $n\times n$ components of $Q$ are the frontal slices $\mathcal{Q}_i \in \R^{n\times n}$ of the tensor $\mathcal{Q}$, i.e. 
$Q = \begin{bmatrix}  \mathcal{Q}_1 & \cdots & \mathcal{Q}_n\end{bmatrix}$. The mode-2 and mode-3 matricization can be defined as
\bes
\begin{aligned}
Q^{(2)} &= \begin{bmatrix}  \mathcal{Q}_1^T & \cdots & \mathcal{Q}_n^T\end{bmatrix},\\
Q^{(3)} &= \begin{bmatrix}  vec(\mathcal{Q}_1) & \cdots & vec(\mathcal{Q}_n)\end{bmatrix}^T.
\end{aligned}
\ees
It is observed that the following property holds
\be \label{2Hproperty}
w^TQ(u\otimes v)=u^TQ^{(2)}(v\otimes w),
\ee
where $w,u,v\in \Rn$ are arbitrary and $Q$ is symmetric in the sense that $Q(u\otimes v)=Q(v\otimes u)$, see \cite{morKolB09}.  
Let $G(s):= sE-A$, then by using 
\bes
\frac{\partial G(s)^{-1} }{\partial s}= -G(s)^{-1}\frac{\partial G(s)}{\partial s}  G(s)^{-1},
\ees
and \eqref{2Hproperty}, we have
\be
\begin{aligned}
\frac{\partial H_2(s_1,s_2) }{\partial s_1}= -y_1&(s_1+s_2)^{T}Ex_2(s_1,s_2)\\
&-x_1(s_1)^TE^Ty_2(s_1,s_2)
\end{aligned}
\ee
where $y_1(s):= (sE-A)^{-T}C^T$ and $y_2(s_1, s_2):= (s_1E-A)^{-T}C(s_1,s_2)^T$ in which 
$$
C(s_1,s_2) = Q^{(2)}\big(x_1(s_2)\otimes y_1(s_1+s_2)\big)+\frac{1}{2!}N^Ty_1(s_1+s_2)
$$
Similarly
\be
\begin{aligned}
\frac{\partial H_2(s_1,s_2) }{\partial s_2}= -y_1&(s_1+s_2)^{T}Ex_2(s_1,s_2)\\
&-x_1(s_2)^TE^Ty_2(s_2,s_1)
\end{aligned}
\ee
Notice that when $s_1=s_2=\sigma$, the two partial differentiations are the same. This condition on interpolation points is assumed in \cite{morBenB15} to show the moment-matching properties of the ROM.
In the following, we show moment-matching in the multivariate settings where $s_1\neq s_2$ ($s_1=\sigma_{1i}$ and $s_2=\sigma_{2i}$).

\subsection{Moment-Matching for QBDAE}
The goal of a moment-matching based reduction approach is to ensure that the high-order transfer functions are well approximated. In case of symmetric transfer functions, we can represent it as
\be \label{goal_tf}
H_k(s_1,\ldots,s_k)\approx \hat{H}_k(s_1,\ldots,s_k), \quad \mbox{for }k=1,\ldots,K,
\ee
with $\hat{H}_k(s_1,\ldots,s_k)$ being the k-th order multivariate transfer function of the reduced system \eqref{reduced_qbdae}. With the task in \eqref{goal_tf} achieved for some $K$, we can expect that the output $y(t)$ is well approximated by $\hy(t)$. To get recursive relations between vectors for approximation subspaces, it is assumed in \cite{morBenB15}  that $s_1=s_2=\sigma$. With these settings, the second-order transfer function becomes
\bes
H_2(\sigma,\sigma)= y(2\sigma)^T\Big( Q \left(x_1(\sigma)\otimes x_1(\sigma)\right) + Nx_1(\sigma)\Big).
\ees 

The following Lemma summarizes the result introduced in \cite{morBenB15}.

\begin{lemma} \label{lemma_sec2}
Let $\sigma_i\in \C$ be the interpolation points and $\sigma_i\notin \{\Lambda(A,E), \Lambda(A_r,E_r)\}$, where $\Lambda(A,E)$ represents the generalized eigenvalues of the matrix pencil $\lambda E-A$. Assume that $\hE=W^TEV$ is nonsingular and $\hA$, $\hQ$, $\hN$, $\hB$, $\hC$ are as in \eqref{reduced_mats} with full rank matrices $V,W\in \Rnr$ such that 
\bes
\begin{aligned}
&\starspan (V)=\starspan_{i=1,\ldots,k}\{x_1(\sigma_i), ~x_2(\sigma_i,\sigma_i)\}, \\
&\starspan(W)=\starspan_{i=1,\ldots,k} \{y_1(2\sigma_i),~y_2(\sigma_i,\sigma_i)\},
\end{aligned}
\ees
then the reduced QBDAE satisfies the following (Hermite) interpolation conditions:
\bes
\begin{aligned}
H_1(\sigma_i)&=\hat{H}_1(\sigma_i), \qquad\quad  H_1(2\sigma_i)=\hat{H}_1(2\sigma_i),\\
H_2(\sigma_i,\sigma_i)&=\hat{H}_2(\sigma_i,\sigma_i), \quad \frac{\partial}{\partial s_j}H_2(\sigma_i,\sigma_i)=\frac{\partial}{\partial s_j}\hat{H}_2(\sigma_i,\sigma_i),~~ j=1,2.
\end{aligned}
\ees
\end{lemma}

See \cite{morBenB15} for a proof. Next, we present moment-matching properties in the multivariable settings, where $s_1\neq s_2$.

\begin{lemma} \label{new_lemma_sec2}
Let $\sigma_{1i},\sigma_{2i}\in \C$ with  $\sigma_{1i},\sigma_{2i}\notin \{\Lambda(A,E), \Lambda(A_r,E_r)\}$. Assume that $\hE=W^TEV$ is nonsingular and $\hA$, $\hQ$, $\hN$, $\hB$, $\hC$ are as in \eqref{reduced_mats} with full rank matrices $V,W\in \Rnr$ such that 
\bes
\begin{aligned}
&\starspan (V)=\starspan_{i=1,\ldots,k}\{x_1(\sigma_{1i}), ~x_1(\sigma_{2i}),~x_2(\sigma_{1i},\sigma_{2i})\}\\
&\starspan(W)=\starspan_{i=1,\ldots,k} \{y_1(\sigma_{1i}+\sigma_{2i}),~y_2(\sigma_{1i},\sigma_{2i}),~y_2(\sigma_{2i},\sigma_{1i})\}.
\end{aligned}
\ees
Then the reduced QBDAE satisfies the following (Hermite) interpolation conditions:
\bes
\begin{aligned}
H_1(&\sigma_{1i})=\hat{H}_1(\sigma_{1i}),   \quad H_1(\sigma_{2i})=\hat{H}_1(\sigma_{2i}),  \quad   H_1(\sigma_{1i}+\sigma_{2i})=\hat{H}_1(\sigma_{1i}+\sigma_{2i}),\\
&H_2(\sigma_{1i},\sigma_{2i})=\hat{H}_2(\sigma_{1i},\sigma_{2i}), \quad \frac{\partial}{\partial s_1}H_2(\sigma_{1i},\sigma_{2i})=\frac{\partial}{\partial s_1}\hat{H}_2(\sigma_{1i},\sigma_{2i}),\\
&\qquad \frac{\partial}{\partial s_2}H_2(\sigma_{2i},\sigma_{1i})=\frac{\partial}{\partial s_2}\hat{H}_2(\sigma_{2i},\sigma_{1i}).
\end{aligned}
\ees
\end{lemma}

The proof of the statement is similar to Lemma~\ref{lemma_sec2} and therefore omitted. Note that the statement in Lemma~\ref{new_lemma_sec2} reduces to Lemma~\ref{lemma_sec2}, if $\sigma_{1i}=\sigma_{2i}$. In the remaining part of the paper, our goal is to identify a good choice of the interpolation points $\sigma_{1i}$ and $\sigma_{2i}$. 

\section{Error Bound for QBDAE's}

In this section, we show how the error bound expression, derived initially in \cite{morFenAB15} for parametric linear time invariant systems, can be extended to the quadratic-bilinear DAEs. We begin with a brief overview of the error bound for the first subsystem, as in \cite{morFenAB15} and then discuss the extension to the second subsystem of QBDAE. 
\subsection{Error bound for $H_1(s_1)$}
Here the error bound provides an estimate for the error between $H_1(s_1)$ and $\hat H_1(s_1)$. To this end, we define the primal and the dual systems as:
\begin{align}\label{primal1}
(s_1E-A)x_1(s_1)&=B,   \\  \label{dual1}  
(s_1E-A)^Tx_1^{du}(s_1)&=-C^T,    
\end{align}
respectively, where $T$ denotes transpose of the matrix. The error bound is constructed so that it is based on two residuals, which result from MOR of the primal and the dual system, respectively. The primal system is reduced using the matrix pair $V_{1}$ and $W_{1}$, where 
\be
\starspan (V_{1})=\starspan_{i=1,\ldots,k}\{x_1(\sigma_{1i})\}, \quad \starspan (W_{1})=\starspan_{i=1,\ldots,k}\{x_1^{du}(\sigma_{1i})\}.   \label{VW1}
\ee
As a result, the reduced primal system is,
\begin{equation*}
(s_1\hE_1-\hA_1)z_1(s_1)=\hB, 
\end{equation*}
where $\hat E_1= W_1^T EV_1$, $\hat A_1= W_1^T AV_1$, $\hat B_1= W_1^T B$ and $\hat C_1= CV_1$. Here $\hx_1(s_1):=V_1z_1(s_1)$ is the approximation of $x_1(s_1)$. Due to the dual relation between \eqref{primal1} and \eqref{dual1}, the dual system can be reduced by using $V_1^{du} = W_{1}$ and $W_{1}^{du} = V_1$. The reduced dual system is
\begin{equation*}
(s_1\tE_1-\tA_1)^Tz_1^{du}(s_1)=-\tC_1^T, 
\end{equation*}
where $\tE_1= V_1^T EW_1$, $\tA_1= V_1^T AW_1$, $\tC_1= W_1^T C^T$. Also $\tx_1^{du}(s_1):=W_1z_1^{du}(s_1)$ is the approximation of $x_1^{du}(s_1)$.
The residuals associated with the reduction of the primal and the dual systems can be written as
\begin{equation}
\begin{aligned}
&r_1^{pr}(s_1)=B-(s_1E-A)V_1z_1(s_1),\\ & r_1^{du}(s_1)=-C^T-(s_1E-A)^TW_1z_1^{du}(s_1).
\end{aligned}
\end{equation}
With these quantities, the following result provides an a posteriori upper bound on the approximation error, $|H_1(s_1)-\hat{H}_1(s_1)|$:
\begin{theorem}\label{feng_main}\cite{morFenAB15}
The upper bound on the approximation of the transfer function $H_1(s_1)=C(s_1E-A)^{-1}B$ can be written as $|H_1(s_1)-\hat{H}_1(s_1)|=\Delta_1(s_1)$, where 
\begin{equation} \label{delta_1}
\Delta_1(s_1):=\frac{\|r_1^{du}(s_1)\|_2\|r_1^{pr}(s_1)\|_2}{\beta_1(s_1)},
\end{equation}
in which $\beta_1(s_1)=\sigma_{\min} (G(s_1))$, where $\sigma_{\min}$ indicates the smallest singular value of $G(s_1)$. 
\end{theorem}

\subsection{Error Bound for $H_2(s_1,s_2)$}
Analogous to $H_1(s_1)$, we define the primal and dual systems as:
\begin{align}\label{primal}
&G(s_{1}+s_2)x_2(s_{1},s_2)=B(s_1,s_2),\\ \label{dual}    
&G^T(s_{1}+s_2)x_2^{du}(s_{1},s_2)=-C^T,  
\end{align}
respectively. The interpolation points for $H_1(s_1)$ can be identified through the error bound $\Delta_1(s_1)$ by using a greedy framework as presented in \cite{morFenAB15}. This means that we can select  $\sigma_{1i}$ for $i=1,\ldots,r$ as the interpolation points corresponding to the maximal values of the error bound at subsequent iterations of the greedy algorithm in  \cite{morFenAB15}.With these interpolation points fixed for $s_1$, we can also express error bound for the second subsystem. The error bound is constructed based on two residuals, which result from MOR of the primal and the dual systems in \eqref{primal} \eqref{dual}, respectively. The primal system is reduced using the matrix pair $V_{2}$ and $W_{2}$, where 
\be
\starspan (V_{2})=\starspan_{j=1,\ldots,k}\{x_2(\sigma_{1i},\sigma_{2j})\}, \quad \starspan (W_{2})=\starspan_{j=1,\ldots,k}\{x_2^{du}(\sigma_{1i},\sigma_{2j})\}. \label{VW2}
\ee
As a result, the reduced primal system is
\begin{equation*}
((s_1+s_2)\hE_2-\hA_2)z_2(s_1,s_2)=\hB(s_1,s_2),
\end{equation*}
where $\hat E_2= W_2^T EV_2$, $\hat A_2= W_2^T AV_2$, $\hat B(s_1,s_2)= W_2^T B(s_1,s_2)$ and $\hat C_2=  CV_2$. Similarly, the dual system is reduced using the matrix pair $V_{2}^{du}$ and $W_{2}^{du}$, 
\be
\starspan (V_{2}^{du})=\starspan_{i=1,\ldots,k}\{x_2^{du}(\sigma_{1i},\sigma_{2i})\}, \quad \starspan (W_{2}^{du})=\starspan_{i=1,\ldots,k}\{x_2(\sigma_{1i},\sigma_{2i})\}. 
\ee
The reduced dual system is
\begin{equation*}
((s_1+s_2)\tE_2-\tA_2)^Tz_2^{du}(s_1,s_2)=-\tC_2^T,
\end{equation*}
where $\tE_2= (W_2^{du})^T EV_2^{du}$, $\tA_2= (W_2^{du})^T AV_2^{du}$, $\tC^T_2= (V_2^{du})^T C^T$.
The residuals associated with the reduction of the primal and dual systems can be written as
\begin{equation}
\begin{aligned}
&r_2^{pr}(s_1,s_2)=B(s_1,s_2)-((s_1+s_2)E-A)V_2z_2(s_1,s_2),\\ & r_2^{du}(s_1,s_2)=-C^T-((s_1+s_2)E-A)^TV_2^{du}z_2^{du}(s_1,s_2).
\end{aligned}
\end{equation}
With these quantities, the following result provides an a posteriori upper bound on the approximation error, $|H_2(s_1,s_2)-\hat{H}_2(s_1,s_2)|$:
\begin{theorem}
The upper bound on the approximation of $H_2(s_1,s_2)=C((s_1+s_2)E-A)^{-1}B(s_1,s_2)$ can be written as $|H_2(s_1,s_2)-\hat{H}_2(s_1,s_2)|=\Delta_2(s_1,s_2)$, where 
\begin{equation}\label{delta_2}
\Delta_2(s_1,s_2):=\frac{\|r_2^{du}(s_1,s_2)\|_2\|r_2^{pr}(s_1,s_2)\|_2}{\beta_2(s_1,s_2)},
\end{equation}
in which $\beta_2(s_1,s_2)=\sigma_{\min} (G(s_1+s_2))$, where $\sigma_{\min}$ indicates the smallest singular value of $G(s_1+s_2)=(s_1+s_2)E-A$. 
\end{theorem}
The proof is similar to Theorem \ref{feng_main} and therefore is omitted. 

\section{Interpolation Points using Error Bounds}
As discussed in Section~2, the projection matrices $V$ and $W$ defined in Lemma~\ref{new_lemma_sec2} require a good choice of interpolation points $\sigma_{1i}$ and $\sigma_{2i}$ which also serve as interpolation points for MOR of the primal and dual systems in \eqref{primal1}-\eqref{dual1} and  \eqref{primal}-\eqref{dual}. In this section, we show the use of the error bound expressions derived previously to select the interpolation points.  


The idea is to identify interpolation points corresponding to the maximal bound $\Delta_1 (s_1)$. Assuming that $\sigma_{1i}$ are the selected interpolation points for $s_1$, the remaining interpolation points for $s_2$ correspond to the maximal bound $\Delta_2(\sigma_{1i},s_2)$ for each value of $\sigma_{1i}$. In this way, the error bound can be used iteratively to select a good choice of interpolation points in a predefined sample space, starting from an initial choice of sigma's. The selected interpolation points are then used to construct and update the required basis matrices $V$ and $W$, by using the multimoment-matching technique described before. It is interesting to see that although we need to construct the ROMs for the primal and the dual systems in \eqref{primal1}-\eqref{dual1} and  \eqref{primal}-\eqref{dual}, the projection matrices for those ROMs are obtained without extra computations, since   $V_1, W_1$ and $V_2, W_2$ are part of $V, W$ by definition. Therefore, $V, W$ can be obtained by orthogonalizing $V_1$ with $V_2$ and $W_1$ with $W_2$ as indicated in Step \ref{orth} of Algorithm~\ref{CHalgorithm}, where a greedy framework for selecting interpolation points is presented. For an initial pair of interpolation points, the ROMs of the primal and the dual systems in  \eqref{primal1}-\eqref{dual1} and  \eqref{primal}-\eqref{dual} are constructed and the error bounds $\Delta_1, \Delta_2$ are computed. A new pair is selected such that the corresponding error bounds $\Delta_1$ and $\Delta_2$ are maximized at these points. With the selected interpolation points, we enrich the projection matrices $V, W$ for MOR of the original quadratic-bilinear system iteratively during the greedy algorithm. Finally, the  reduced quadratic bilinear system is constructed using $V, W$ that are derived upon convergence of Algorithm~\ref{CHalgorithm}.  Algorithm~\ref{CHalgorithm} stops when $\Delta:=\Delta_1+\Delta_2$ is below the tolerance $\epsilon_{tol}$, where $\Delta$ includes the errors introduced by approximating the first and second transfer functions. Since the interpolation points are selected according to the error bounds $\Delta_1$ and $\Delta_2$, it is important that the error bounds dynamically reflect the decay of the true error with the iteration of the greedy algorithm. Ideally, the error bounds should be very close to the true error. Numerical tests in the next section show that the error bounds really control the true error robustly.

    \scalebox{0.95}{
\begin{minipage}{1\linewidth}
\begin{algorithm}[H]
\caption{An adaptive framework for selection of interpolation points}
\label{CHalgorithm}
{\bf Inputs}:  $\sigma_{10}$, $\sigma_{20}$, $E$, $A$, $N$, $H$, $B$, $C$ and $S_{{\textrm{sample}}}$: a set of the samples of $\mu:=(s_1,s_2)$, which covers the domain of the two frequency variables. \\
{\bf Outputs}: $\mu$, $V$ and $W$ \\
Initialization: $V=[~]$;  $W=[~]$; $V_1=[~]$; $W_1=[~]$; $V_2=[~]$; $W_2=[~]$; $\epsilon=1$; $i=-1$; $j=0$; $\epsilon_{tol} <1$, $\mu^0=(\sigma_{10}, \sigma_{20})$.  \\
WHILE $\epsilon>\epsilon_{tol}$
\begin{enumerate}[label=\arabic*]
\item \quad $i=i+1$;  $j=j+1$;
\item \quad compute $\bar{V}_i({\sigma_{1i}})$ and $\bar{W}_i (\sigma_{1i})$ using \eqref{VW1}   
\item  \quad $V_1=orth[V_1, \bar{V}_i]$;~ $W_1=orth[W_1,\bar{ W}_{i}]$;
\item \quad $\sigma_{1j}=\arg \max\limits_{\sigma_1 \in S_{1\  \textrm{sample}}} \Delta_1(\sigma_1)$;
\item \quad compute $V_i({\sigma_{1i}},{\sigma_{2i}})$ and $W_i (\sigma_{1i},{\sigma_{2i}})$  using  \eqref{VW2}
\item  \quad $V_2=orth[V_2, {V}_i]$;~ $W_2=orth[W_2,{W}_{i}]$; 
\item \quad $\sigma_{2j}=\arg \max\limits_{\sigma_2 \in S_{2\  \textrm{sample}}} \Delta_2(\sigma_{1i},\sigma_2)$;
\item \quad $\mu^j = [\sigma_{1j},\sigma_{2j}]$;
\item \label{orth} \quad $V=orth[V_1, V_2 ]$;~ $W=orth[W_1,W_2]$;
\item \quad $\Delta(\mu^j) :=  \Delta_1(\mu^j) + \Delta_2(\mu^j)$; ~ ${\epsilon}=\Delta(\mu^j)$;

END WHILE.
\end{enumerate}
\end{algorithm}
\end{minipage}%
    }

\section{Numerical results}\label{section_results}
We consider three benchmark examples for our results on MOR of QBDAE systems. The results are compared with the one-sided and two-sided projection methods, where the interpolation points are computed by IRKA, implemented on the linear part of the system. We represented the proposed method by 1s/2s-greedy(One-sided/two-sided projection with greedy based interpolation points) and the method from literature by 1s/2s-IRKA (One-sided/two- sided projection with IRKA interpolation points). The Max. True error in the following tables is defined as $\max\limits_{s_1, s_2 \in  S_{\textrm{sample}}} |H_1(s_1)-\hat H_1(s_1)|+|H_2(s_1,s_2)-\hat H_2(s_1,s_2)|$ and the Max. error bound is $\max\limits_{s_1, s_2 \in S_{\textrm{sample}}}  \Delta(s_1, s_2)$.
\subsection{Nonlinear RC circuit}   
\begin{figure}[H]
\centering
\resizebox{9cm}{!}{\begin{tikzpicture}

\draw (0cm,3cm) -- (3cm,3cm);

\node at (1.2cm,3cm) [above] {{\small $v_1$}};

\node[rectangle,draw=black,minimum width=0.6cm,
    minimum height=.3cm,inner sep=0] at (3.3cm,3cm) {};
\node at (3.3cm,3cm) [above,yshift=3] {{\small $g(v)$}};
\draw (3.6cm,3cm) -- (4.2cm,3cm);

\node at (4.2cm,3cm) [above] {{\small $v_2$}};
\draw (4.2cm,3cm) -- (6.2cm,3cm) [dashed];

\node at (6.2cm,3cm) [above] {{\small $v_{l-2}$}};
\draw (6.2cm,3cm) -- (6.8cm,3cm);

\node[rectangle,draw=black,minimum width=0.6cm,
    minimum height=.3cm,inner sep=0] at (7.1cm,3cm) {};
\node at (7.1cm,3cm) [above,yshift=3] {{\small $g(v)$}};

\node at (8cm,3cm) [above] {{\small $v_{l-1}$}};

\draw (7.4cm,3cm) -- (8.6cm,3cm);
\node[rectangle,draw=black,minimum width=0.6cm,
    minimum height=.3cm,inner sep=0] at (8.9cm,3cm) {};
\node at (8.9cm,3cm) [above,yshift=3] {{\small $g(v)$}};
\draw (9.2cm,3cm) -- (9.8cm,3cm);

\node at (9.8cm,3cm) [above] {{\small $v_l$}};

\draw (9.8cm,3cm) -- (9.8cm,1.6cm);
\draw (9.56cm,1.6cm) -- (10.04cm,1.6cm) [thick];
\draw (9.56cm,1.4cm) -- (10.04cm,1.4cm) [thick];
\node at (9.6cm,1.6cm) [left] {{\small $C$}};
\draw (9.8cm,1.4cm) -- (9.8cm,0cm);

\draw (9.8cm,0cm) -- (6.2cm,0cm);
\draw (6.2cm,0cm) -- (4.2cm,0cm) [dashed];
\draw (4.2cm,0cm) -- (0cm,0cm);

\draw (0cm,0cm) -- (0cm,1.25cm);
\node[circle,draw=black,minimum width=0.5cm] at (0,1.5cm) {};
\draw[->] (0,1.32cm) --(0,1.66cm){};
\node at (0,1.5cm) [left,xshift=-7,yshift=-11]
{\begin{rotate}{90}{\small $i=u(t)$}\end{rotate}};
\draw (0cm,1.75cm) -- (0cm,3cm);

\draw (1.2cm,0cm) -- (1.2cm,1.2cm);
\node[rectangle,draw=black,minimum width=.3cm,
    minimum height=0.6cm,inner sep=0] at (1.2cm,1.5cm) {};
\node at (1.2cm,1.5cm) [left,xshift=-5,yshift=-7]
{\begin{rotate}{90}{\small $g(v)$}\end{rotate}};
\draw (1.2cm,1.8cm) -- (1.2cm,3cm);

\draw (2.4cm,0cm) -- (2.4cm,1.4cm);
\draw (2.16cm,1.4cm) -- (2.64cm,1.4cm) [thick];
\draw (2.16cm,1.6cm) -- (2.64cm,1.6cm) [thick];
\node at (2.2cm,1.6cm) [left] {{\small $C$}};
\draw (2.4cm,1.6cm) -- (2.4cm,3cm);

\draw (4.2cm,0cm) -- (4.2cm,1.4cm);
\draw (3.96cm,1.4cm) -- (4.44cm,1.4cm) [thick];
\draw (3.96cm,1.6cm) -- (4.44cm,1.6cm) [thick];
\node at (4cm,1.6cm) [left] {{\small $C$}};
\draw (4.2cm,1.6cm) -- (4.2cm,3cm);

\draw (6.2cm,0cm) -- (6.2cm,1.4cm);
\draw (5.96cm,1.4cm) -- (6.44cm,1.4cm) [thick];
\draw (5.96cm,1.6cm) -- (6.44cm,1.6cm) [thick];
\node at (6cm,1.6cm) [left] {{\small $C$}};
\draw (6.2cm,1.6cm) -- (6.2cm,3cm);

\draw (8cm,0cm) -- (8cm,1.4cm);
\draw (7.76cm,1.4cm) -- (8.24cm,1.4cm) [thick];
\draw (7.76cm,1.6cm) -- (8.24cm,1.6cm) [thick];
\node at (7.8cm,1.6cm) [left] {{\small $C$}};
\draw (8cm,1.6cm) -- (8cm,3cm);

\draw (9.8cm,0cm) -- (10.2cm,0cm);
\draw (10.2cm,0cm) -- (10.2cm,-0.5cm);
\draw (9.8cm,-.5cm) -- (10.6cm,-.5cm);
\draw (10cm,-.6cm) -- (10.4cm,-.6cm);


    minimum height=0.6cm,inner sep=0] at (0.9375cm,1.3125cm) {};





    minimum height=.225cm,inner sep=0] at (6.09375cm,2.625cm) {};

    minimum height=.225cm,inner sep=0] at (7.03125cm,2.625cm) {};


    minimum height=0.5175cm,inner sep=0] at (7.5cm,1.35cm) {};



\end{tikzpicture}}
\caption{Nonlinear RC circuit}
 \label{fig:RC_circuit}
\end{figure}
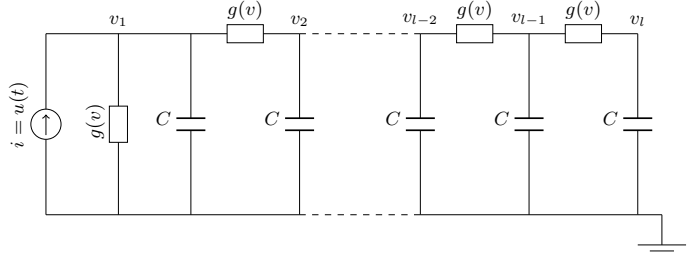

The nonlinear RC circuit was first considered in \cite{morChe99} and since then it has been used in many papers for nonlinear  MOR \cite{morFre11}. Consider $v$ be the voltage and $g(v)$ be the current function then I-V characteristics can be represented as:  $g(v)= e^{40v} + v -1$.  The nonlinearity in the current function results in nonlinear model. All the capacitances are fixed to $C=1$. Figure~\ref{fig:RC_circuit} shows the complete circuit.

It is shown in \cite{morGu11} that the nonlinearity in the RC circuit can be written in the quadratic-bilinear form as in \eqref{original_qbdae} by introducing some auxiliary variables. The transformation is exact, but the dimension of the system increases to $n =2\cdot l$, where $l$ represents the number of nodes in Figure~\ref{fig:RC_circuit}, and it is also the dimension of the original nonlinear system.  

For our results, we set $l=50$, so $n=100$ and use two-sided projection method to reduce the system.
~Table~\ref{table1:RCq} shows the results with tolerance $\epsilon_{tol}=1e^{-5}$  and an initial choice of interpolation points as $\sigma_{1}=\sigma_{20}=119.5642$.

\begin{table}[H]
\begin{center}
\begin{tabular}{ c c c c }
S.No.\quad & Interpolation points \quad& Max. True Error  &  Max. Est. Error\\
& $\{\sigma_{1i},\sigma_{2i}\}$& &\\
1&$119.5642,119.5642$&$1.8616\times10^{-2}$&$0.109183$\\
2&$0.9875,0.9875$&$1.3683\times10^{-3}$&$8.4421\times10^{-3}$\\
3&$4.9567,0.9875$&$1.6127\times10^{-4}$&$4.0341\times10^{-4}$\\
4&$18.1107,5.5319$&$4.2956\times10^{-5}$&$7.22\times10^{-5}$\\
5&$2.0292,4.4445$&$8.239\times10^{-6}$&$9.6404\times10^{-6}$\\
\end{tabular}
\end{center}
\caption{Error estimation results for RC circuit}
\label{table1:RCq}
\end{table}

\begin{figure}[H]
\begin{subfigure}{.5\textwidth}
  \centering
  \includegraphics[width=.8\linewidth]{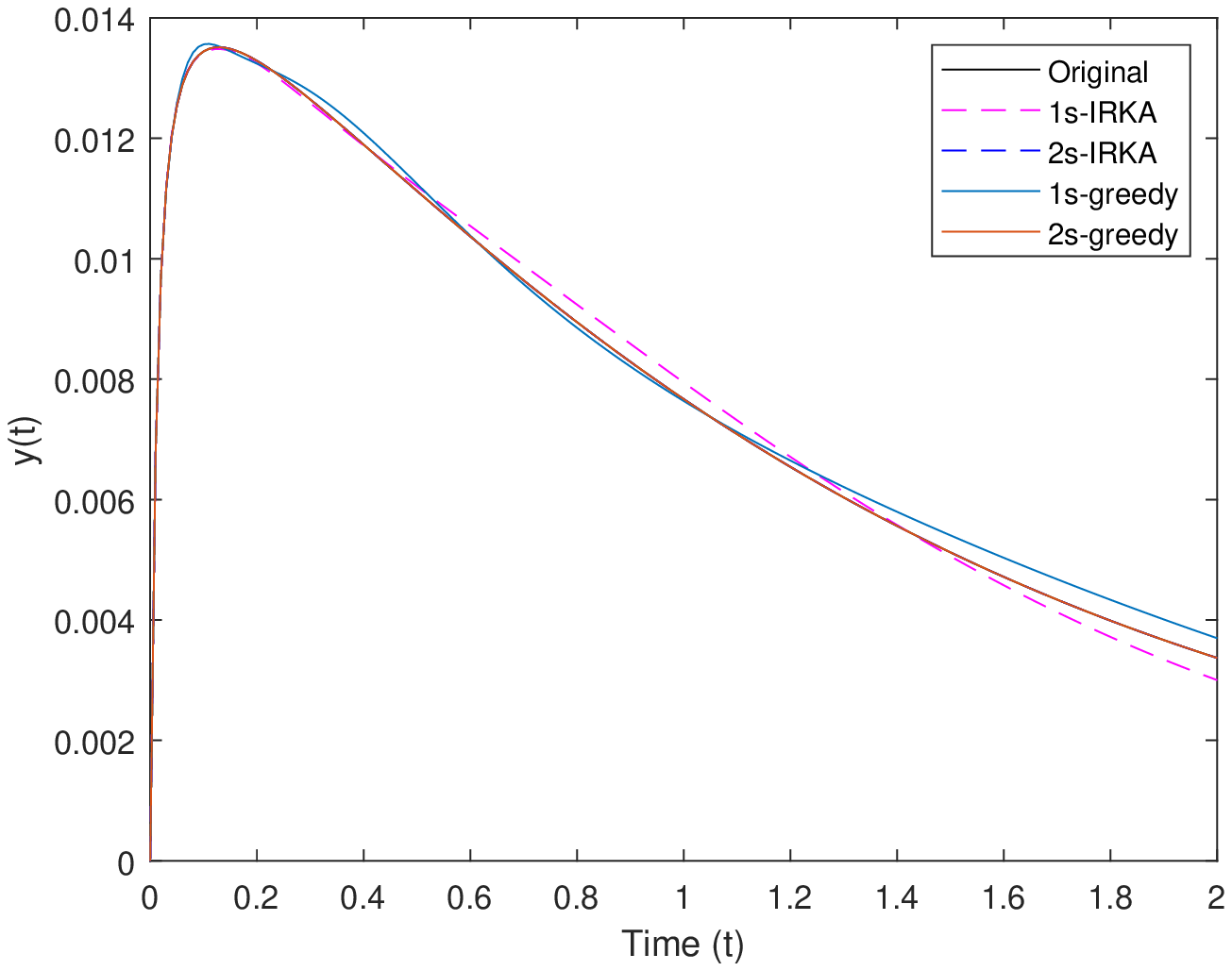}  
  \caption{Comparison of transient response}
  \label{fig:Transient_RC}
\end{subfigure}
\begin{subfigure}{.5\textwidth}
  \centering
  \includegraphics[width=.8\linewidth]{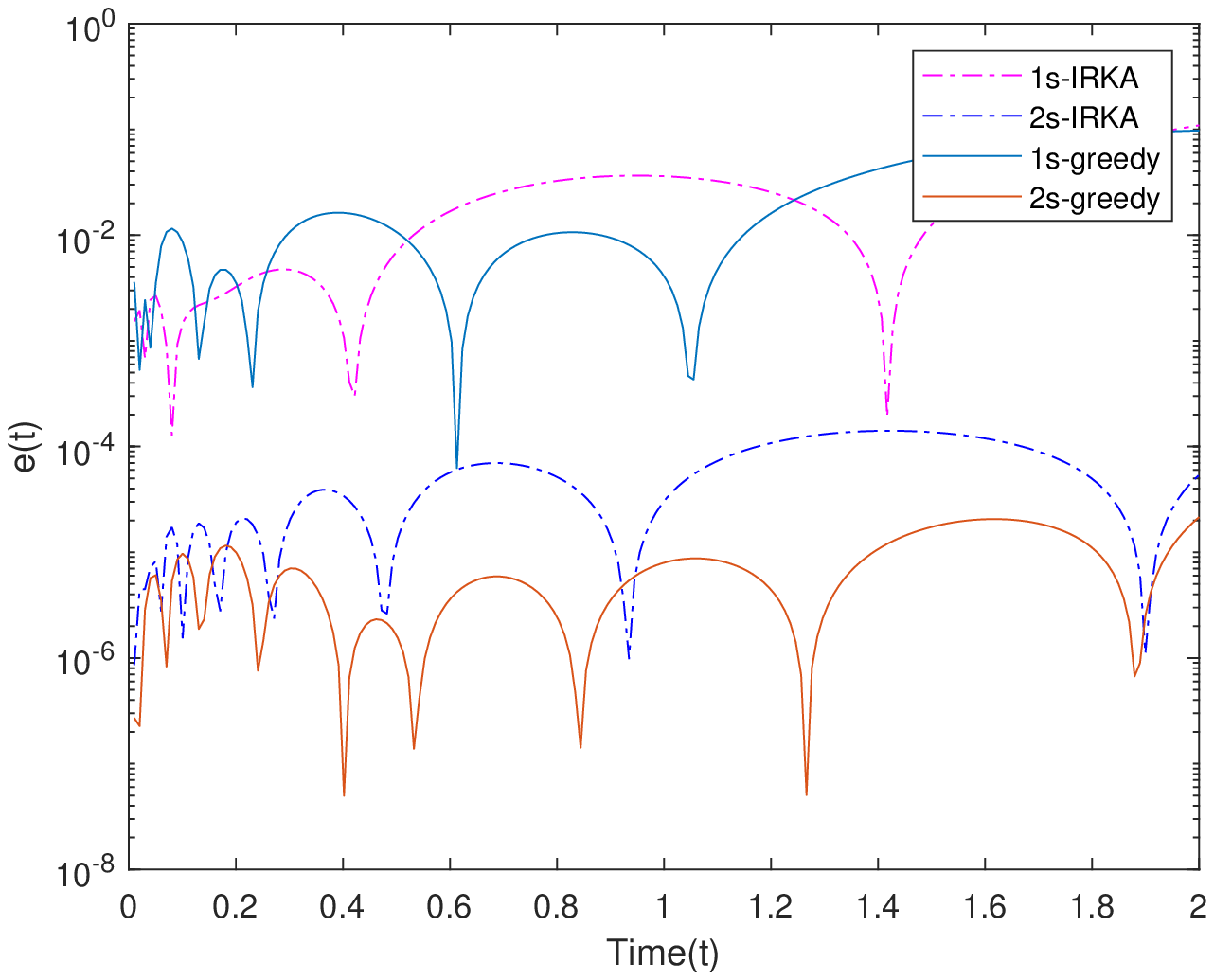}  
  \caption{Comparison of relative error}
  \label{fig:Rel_Error_RC}
\end{subfigure}
\caption{Non-linear RC circuit}
\label{fig:RC_Circuit}
\end{figure}

The second column of Table~ \ref{table1:RCq} shows interpolation points that are identified by the greedy framework and are based on the error bound. It is clear that the error bound tightly catches the true error and can be used as a surrogate of the true error to select the interpolation points. The size of the ROM obtained from both approaches has been kept the same i.e. $r_1 = r_2 = 12$. For the input $u(t) = e^{-t}$, the output of the original model and ROMs along with corresponding relative errors are shown in~Figure \ref{fig:RC_Circuit}.

Figure~\ref{fig:Transient_RC} shows the comparison of transient response of the two approaches, while Figure~\ref{fig:Rel_Error_RC} plots relative errors of the two approaches. It is clearly seen that 1s-greedy and 2s-greedy outperform 1s-IRKA and 2s-IRKA respectively in terms of accuracy.

\subsection{Burgers' Equation}

In nonlinear MOR, 1D  burgers' equation is commonly used \cite{morKunV08},\cite{morBenB15}. Mathematical model of 1D burger's equation  with $\Gamma = (0,1) \times (1,T)$  is:
\begin{equation}
\begin{aligned}
\upsilon_t + \upsilon\upsilon_x = \nu \cdot \upsilon\upsilon_{xx},\ \ \ \ \ \ \ \ \ \ \ \ \ \ \ \ \  in\  \Gamma, \\
\alpha\upsilon(0,t) +\beta x(0,t) = u(t), \ \ \upsilon_x(1,t) = 0, \ \ \ \ \ \ \ \ \ \  t\in (0,T), \\
\upsilon(x,0) = \upsilon_0(x), \ \ \upsilon_0(x) = 0,  \ \ \ \ \ \ \ \ \  x\in (0,1), \\
\end{aligned}
\label{burger_eq}
\end{equation}
we use it as an example to test our proposed method. We keep the size of the original model as n = 1000.~Table~\ref{table2:Burg}  shows our results with tolerance $\epsilon_{tol}=1e^{-4}$ and an initial choice of interpolation points as $\sigma_{10} = \sigma_{20} = 5.4124$.

\begin{table}[H]\label{Burg_table2}
\begin{center}
\begin{tabular}{ c c c c }
S.No.\quad & Interpolation points \quad& Max. True Error  &  Max. Est. Error\\
& $\{ \sigma_{1i}, \sigma_{2i}\}$& &\\
1&$5.4124 ,5.4124 $&$1.1299\times10^{-3}$&$32.4786$\\
2&$31.6141,1.383 $&$1.0259\times10^{-3}$&$3.2407$\\
3&$2.9603 ,1.0818 $&$1.0746\times10^{-3}$&$4.2125\times10^{-1}$\\
4&$9.2633 -11.3351\iota,24.9534$&$1.416\times10^{-4}$&$4.3411\times10^{-4}$\\
5&$7.4119 - 3.622\iota,1.0818$&$1.785\times10^{-5}$&$1.7869\times10^{-5}$\\
\end{tabular}
\end{center}
\caption{Error estimation results for burgers equation}
\label{table2:Burg}
\end{table}

\begin{figure}[H]
\begin{subfigure}{.5\textwidth}
  \centering
  \includegraphics[width=.8\linewidth]{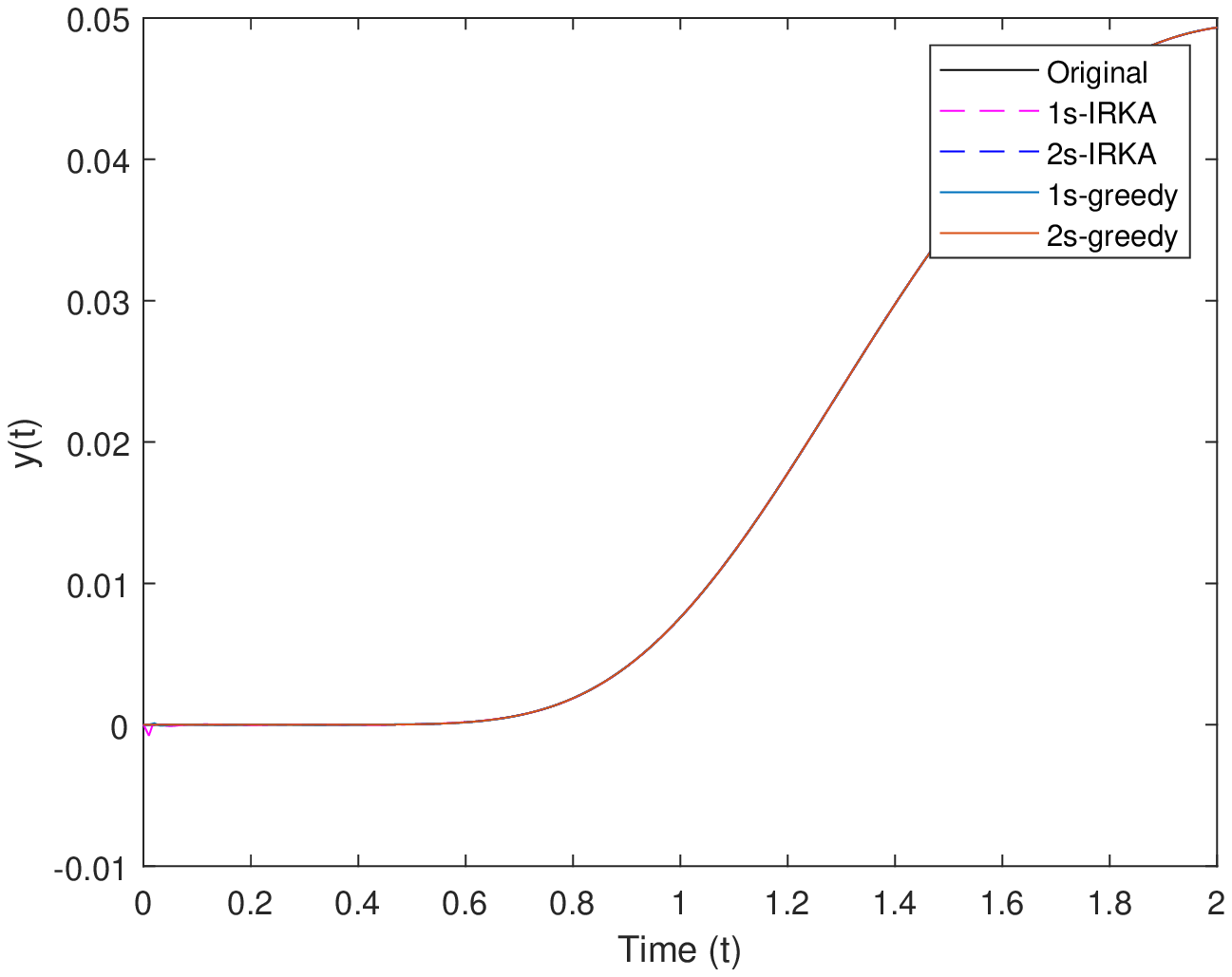}  
  \caption{Comparison of transient response}
  \label{fig:Transient_Burg}
\end{subfigure}
\begin{subfigure}{.5\textwidth}
  \centering
  \includegraphics[width=.8\linewidth]{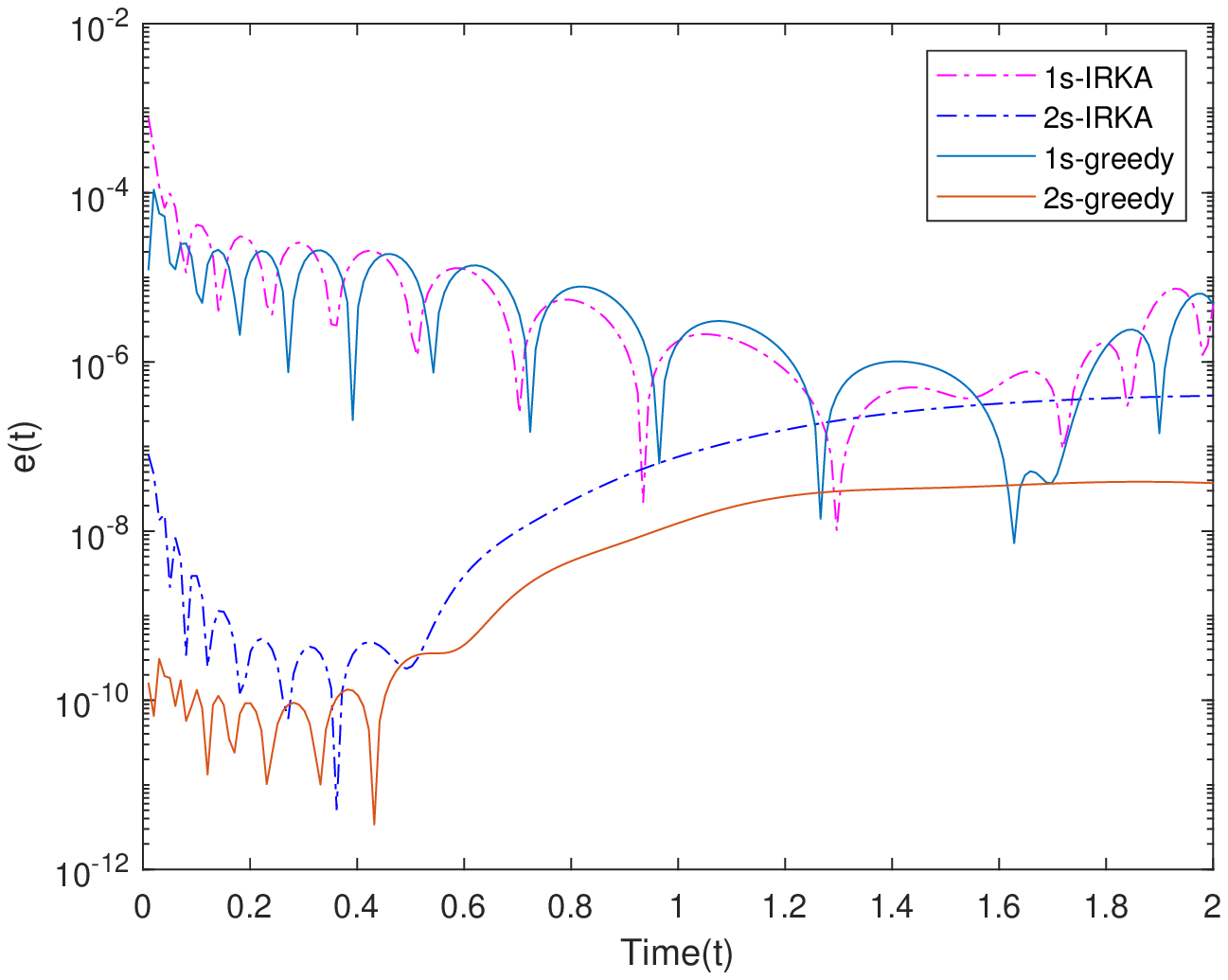}  
  \caption{Comparison of absolute error}
  \label{fig:Abs_Error_Burg}
\end{subfigure}
\caption{Burger's equation}
\label{fig:Burg}
\end{figure}

The second column of the table shows interpolation points that are based on the error bound and identified by the greedy framework. Similarly, the error bound again tightly bounds the true error and therefore is reliable for choosing the interpolation points in the greedy algorithm. The  sizes of the ROMs obtained from both approaches are kept same i.e. $r_1 = r_2 = 16$. The ROMs constructed from IRKA interpolation points and the proposed framework are shown in Figure \ref{fig:Burg} for input $u(t)  =  cos(\pi t)$.

Figure~\ref{fig:Transient_Burg} shows the transient responses of the burgers equation computed from simulating the original model and the two different MOR approaches, while Figure \ref{fig:Abs_Error_Burg} compares  the absolute response errors of the ROMs derived using two approaches. The absolute error of ROM constructed using the proposed methodology of choosing interpolation points is less than that of the ROM constructed using IRKA interpolation points, especially for the two-sided projection. 

\subsection{FitzHugh - Nagumo System}

We use the FitzHugh - Nagumo system as our third example to check our results. The FitzHugh - Nagumo system can be represented as\cite{morBenGG16}:

\begin{equation}
\begin{aligned}
\epsilon\upsilon_t(x,t) = \epsilon^{2}\upsilon_{xx}(x,t) + f(\upsilon(x,t)) - w(x,t) + g , \\
w_t(x,t) = h\upsilon(x,t) -\gamma w(x,t) +g, \\
\end{aligned}
\label{FitzH_eq}
\end{equation}

with\ $f(\upsilon) = \upsilon(\upsilon -0.1)(1-\upsilon)$\ and\ boundary\ conditions:

\begin{equation}
\begin{aligned}
\upsilon(x,0) = 0, \ \ \ \ \ \ \ \ \ \ \ \ w(x,0) = 0, \\
\upsilon_x(0,t) = - i_0(t), \ \ \ \ \ \ \ \ \ \upsilon_x(1,t) = 0.\\
\end{aligned}
\end{equation}

Here, we choose $\epsilon = 0.015$, $h=0.5$, $\gamma = 0.05$ and $i_0(t) = 5 \times10^4 t^3 e^{-15t}$. When standard finite difference method is applied to numerically discretize the PDEs in \eqref{FitzH_eq}, a system of  ODEs with cubic non-linearities is obtained. We can get a quadratic bilinear system by introducing new variables. For an original discretized system with size $\bar{n}$, a quadratic bilinear system has the size of $n = 3\bar{n}$. we set $\bar{n} = 100$, which gives rise to quadratic bilinear system of  $n = 300$. We choose interpolation points using the proposed greedy framework to construct the ROM of size $r = 26$ and then compare it with the ROM of the same size, which is constructed from the interpolation points using IRKA.~Table~\ref{table3:FitzH} shows our results with tolerance $\epsilon_{tol}=1e^{-6}$ and the interpolation points $\sigma_{10}=\sigma_{20}=534.69$.

\begin{table}[H]\label{FitzH1}
\begin{center}
\begin{tabular}{ c c c c }
S.No.\quad & Interpolation points \quad& Max. True Error  &  Max. Est. Error\\
& $\{\sigma_{1i},\sigma_{2i}\}$& &\\
1&$534.69,534.69$&$0.282519 $&$1152.4511$\\
2&$1.38,1.08$&$4.7413\times10^{-1}$&$8.4587$\\
3&$3.91-5.45\iota,1.38$&$1.2373\times10^{-4}$&$4.3284 \times10^{-3}$\\
4&$39.38,1.08$&$2.5379 \times10^{-6}$&$5.9555 \times10^{-5}$\\
5&$110.46,1.08$&$8.2393 \times10^{-6}$&$2.1293\times10^{-5}$\\
6&$3.96,1.08$&$4.3429 \times10^{-5}$&$7.1251\times10^{-4}$\\
7&$17.63,1.08$&$7.6047 \times10^{-6}$&$4.6707\times10^{-5}$\\
8&$4.83-4.72\iota,1.08$&$9.7775 \times10^{-8}$&$1.932\times10^{-7}$\\
\end{tabular}
\end{center}
\caption{Error estimation results for the FitzHugh - Nagumo model}
\label{table3:FitzH}
\end{table}

\begin{figure}[H]
\begin{subfigure}{.5\textwidth}
  \centering
  \includegraphics[width=.8\linewidth]{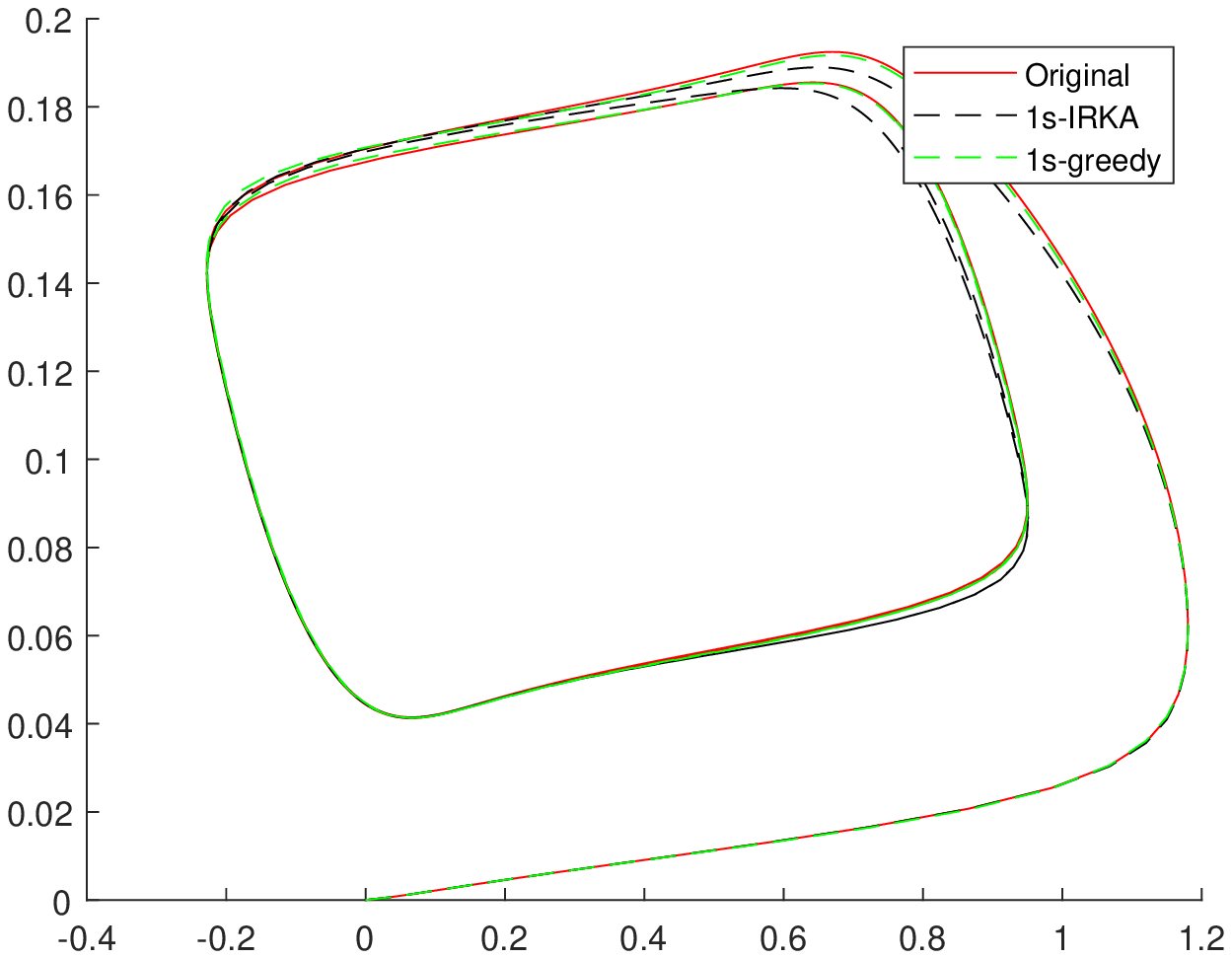}  
  \caption{Comparison of transient response}
  \label{fig:Transient_Fitz}
\end{subfigure}
\begin{subfigure}{.5\textwidth}
  \centering
  \includegraphics[width=.8\linewidth]{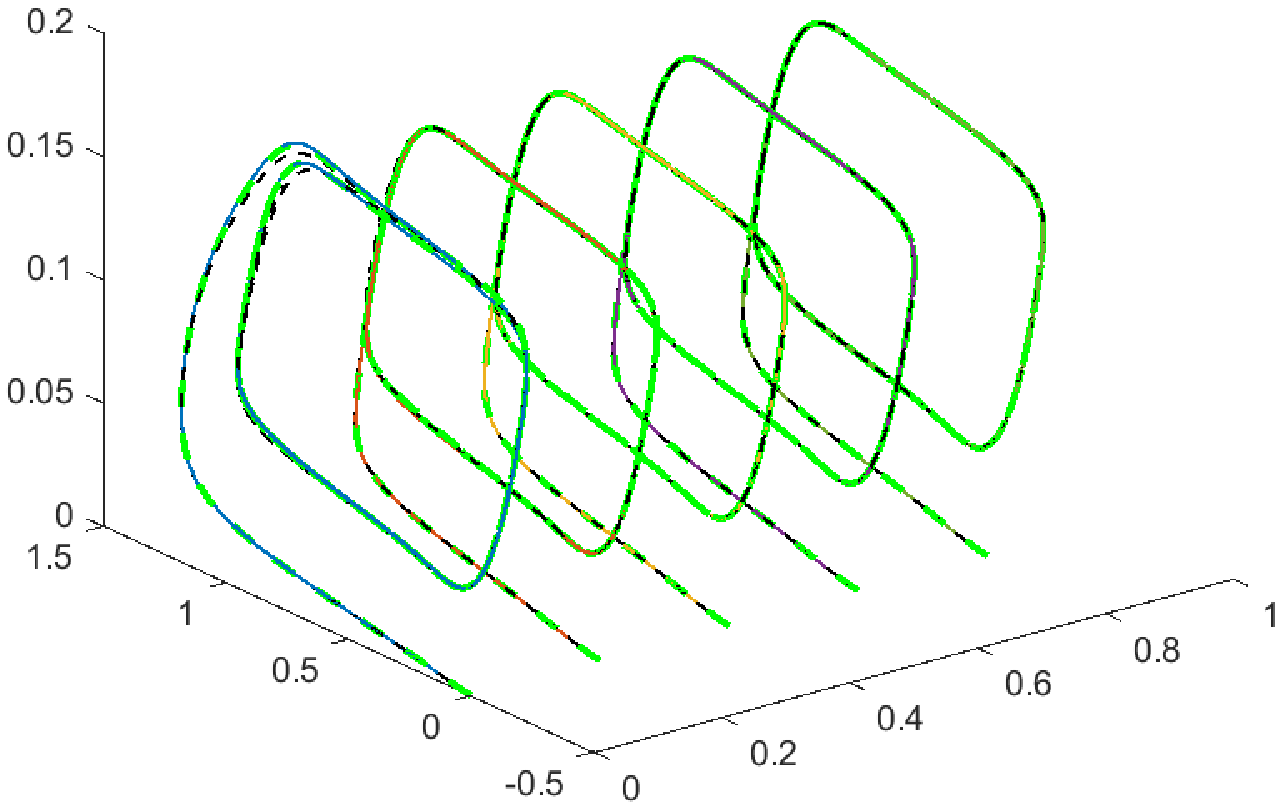}  
  \caption{Comparison of transient response (3D)}
  \label{fig:3D_Fitz}
\end{subfigure}
\caption{FitzHugh - Nagumo equation}
\label{fig:Fitz}
\end{figure}

The table ~\ref{table3:FitzH} shows interpolation points that are selected by the error bound and the decay of the true error and the error bound at each iteration of the greedy algorithm. The error bound once more, estimates the true error accurately, implicating that the selected interpolation points indeed nearly corresponds to the largest error. The  sizes of ROMs obtained from both approaches have been kept the same i.e. $r_1 = r_2 = 26$. Figure ~\ref{fig:Fitz} shows the transient responses of the FitzHugh - Nagumo system computed from simulating the original model and two approaches.

The input signal is $u(t) =50000t^3 e^{-15t}$. It is seen that the 1s-greedy performs better than the 1s-IRKA when the outputs in both cases are compared with that of the original model; however, 2s-greedy and 2s-IRKA produce unstable responses.

\section{Conclusions}\label{section:Conclusions}

In this paper, the proposed methodology of choosing interpolation points for construction of ROM of the first- and second-order transfer functions of quadratic-bilinear systems has been tested for three different models. The results have also been compared with ROMs of the same size constructed using the interpolation points chosen by linear IRKA. In each case, the ROMs constructed using interpolation points from the greedy framework yield better approximation of the output than the  ROMs constructed from IRKA.

\newpage

\bibliography{mybibfile}

\end{document}